
\input graphicx.tex
\input xy
\xyoption{all}


\magnification\magstephalf

\voffset0truecm
\hoffset=0truecm
\vsize=23truecm
\hsize=15.8truecm
\topskip=1truecm

\binoppenalty=10000
\relpenalty=10000

\font\tenbb=msbm10		\font\sevenbb=msbm7		\font\fivebb=msbm5
\font\tensc=cmcsc10		\font\sevensc=cmcsc7 	\font\fivesc=cmcsc5
\font\tensf=cmss10		\font\sevensf=cmss7		\font\fivesf=cmss5
\font\tenfr=eufm10		\font\sevenfr=eufm7		\font\fivefr=eufm5


\newfam\bbfam	\newfam\scfam	\newfam\frfam	\newfam\sffam

\textfont\bbfam=\tenbb
\scriptfont\bbfam=\sevenbb
\scriptscriptfont\bbfam=\fivebb

\textfont\scfam=\tensc
\scriptfont\scfam=\sevensc
\scriptscriptfont\scfam=\fivesc

\textfont\frfam=\tenfr
\scriptfont\frfam=\sevenfr
\scriptscriptfont\frfam=\fivefr

\textfont\sffam=\tensf
\scriptfont\sffam=\sevensf
\scriptscriptfont\sffam=\fivesf

\def\bb{\fam\bbfam \tenbb} 
\def\sc{\fam\scfam \tensc} 
\def\sf{\fam\sffam \tensf} 


\font\sezfont=cmbx10 scaled \magstep1
\font\subsectfont=cmbx10 scaled \magstephalf
\font\titfont=cmbx10 scaled \magstep2
\font\autfont=cmcsc10
\font\intfont=cmss10 

\let\no=\noindent
\let\bi=\bigskip
\let\me=\medskip
\let\sm=\smallskip
\let\ce=\centerline

\let\io=\infty
\def\qqquad{\quad\qquad}


\newcount\sectno\sectno=0
\newcount\subsectno\subsectno=0
\newcount\thmno\thmno=0
\newcount\tagno\tagno=0
\newcount\notitolo\notitolo=0
\newcount\defno\defno=0

\def\sect#1\par{
	\global\advance\sectno by 1 \global\subsectno=0\global\defno=0\global\thmno=0
	\vbox{\vskip.75truecm\advance\hsize by 1mm
	\hbox{\centerline{\sezfont \the\sectno.~~#1}}
	\vskip.25truecm}\nobreak}

\def\subsect#1\par{
	\global\advance\subsectno by 1
	\vbox{\vskip.75truecm\advance\hsize by 1mm
	\line{\subsectfont \the\sectno.\the\subsectno~~#1\hfill}
	\vskip.25truecm}\nobreak}
	
\def\appendix#1\par{
	\global\defno=0\global\thmno=0
	\vbox{\vskip.75truecm\advance\hsize by 1mm
	\hbox{\centerline{\sezfont Appendix:~#1}}
	\vskip.25truecm}\nobreak}

\def\defin#1{\global\advance\defno by 1
	\global\expandafter\edef\csname+#1\endcsname%
    {\number\sectno.\number\defno}
    \no{\bf Definition~\the\sectno.\the\defno.}}

\def\thm#1#2{
	\global\advance\thmno by 1
	\global\expandafter\edef\csname+#1\endcsname%
	{\number\sectno.\number\thmno}
	\no{\bf #2~\the\sectno.\the\thmno.}}

\def\Tag#1{\global\advance\tagno by 1 {(\the\tagno)}
    \global\expandafter\edef\csname+#1\endcsname%
    		{(\number\tagno)}}
\def\tag#1{\eqno\Tag{#1}}

\def\rf#1{\csname+#1\endcsname\relax}

\def\proof{\no{\sl Proof.}\enskip}
\def\qedn{\thinspace\null\nobreak\hfill\hbox{\vbox{\kern-.2pt\hrule height.2pt
        depth.2pt\kern-.2pt\kern-.2pt \hbox to2.5mm{\kern-.2pt\vrule
        width.4pt \kern-.2pt\raise2.5mm\vbox to.2pt{}\lower0pt\vtop
        to.2pt{}\hfil\kern-.2pt \vrule
        width.4pt\kern-.2pt}\kern-.2pt\kern-.2pt\hrule height.2pt
        depth.2pt \kern-.2pt}}\par\medbreak}
    \def\qed{\hfill\qedn}
    
\newif\ifpage\pagefalse
\newif\ifcen\centrue

\headline={
\ifcen\hfil\else
\ifodd\pageno
\global\hoffset=0.5truecm
\else
\global\hoffset=-0.4truecm
\fi\hfil
\fi}

\footline={
	\ifpage
		\hfill\rm\folio\hfill
	\else
		\global\pagetrue\hfill
\fi}

\lccode`\'=`\'

\def\bib#1{\me\item{[#1]\enskip}}


\def\C{{\bb C}} 
\def\R{{\bb R}} 
\def\Z{{\bb Z}} \def\N{{\bb N}}

\let\eps=\varepsilon \let\phe=\varphi

\mathchardef\void="083F

\def\diag{\mathop{\rm Diag}\nolimits}

\def\res{{\rm Res}}

\def\invlim{\mathop{\vtop{\offinterlineskip
\hbox{\rm lim}\kern1pt\hbox{\kern-1.5pt$\longleftarrow$}\kern-3pt}
}\limits}
\def\neweq#1$${\xdef #1{(\the\capno.\the\tagno)}
	\eqno #1$$
	\iffinal\else\rsimb#1\fi
	\global \advance \tagno by 1}
\def\neweqa#1$${\xdef #1{{\rm(\the\capno C.\the\tagno)}}
	\eqno #1$$
	\iffinal\else\rsimb#1\fi
	\global \advance \tagno by 1}
\def\newforclose#1{%
	\global\advance\tagno by 1 
    \global\expandafter\edef\csname+#1\endcsname%
    		{(\number\capno.\number\tagno)}
		\hfil\llap{$(\the\capno.\the\tagno)$}\hfilneg}
\def\newforclosea#1{
	\xdef #1{{\rm(\the\capno C.\the\tagno)}}
	\hfil\llap{$#1$}\hfilneg
	\global \advance \tagno by 1
	\iffinal\else\rsimb#1\fi}
\def\forevery#1#2$${\displaylines{\let\neweqa=\newforclosea
	\let\tag=\newforclose\hfilneg\rlap{$\qqquad\forall#1$}\hfil#2\cr}$$}



\ce{\titfont Holomorphic linearization of commuting}

\sm\ce{\titfont germs of holomorphic maps} 


\me\ce{\autfont Jasmin Raissy}
\sm\ce{\intfont Dipartimento di Matematica e Applicazioni}
\ce{\intfont Universit\`a degli Studi di Milano Bicocca}

\ce{\intfont Via Roberto Cozzi 53, 20125 Milano}

\sm\ce{\intfont E-mail: {\tt jasmin.raissy@unimib.it}}
\bi

{\narrower

{\sc Abstract.} Let $f_1, \dots, f_h$ be $h\ge 2$ germs of biholomorphisms of $\C^n$ fixing the origin. We investigate the shape a (formal) simultaneous linearization of the given germs can have, and we prove that if $f_1,\ldots,f_h$ commute and their linear parts are almost simultaneously Jordanizable then they are simultaneously formally
linearizable.  We next introduce a simultaneous Brjuno-type condition and prove that, in case the linear terms of the germs are diagonalizable, if the germs commute and our Brjuno-type condition holds, then they are holomorphically simultaneously linerizable. This answers to a multi-dimensional version of a problem raised by Moser.

}

\bi
\footnote{}{{\sl\hskip-20pt\noindent Mathematics Subject Classification (2010).} Primary 37F50; Secondary 32H50. \hfill \break {\sl Key words and phrases.} Simultaneous linearization problem, resonances, Brjuno condition, commuting germs, discrete local holomorphic dynamical systems.}
\footnote{}{\sm{\hskip-20pt\noindent\sf Partially supported by FSE, Regione Lombardia,  and by the PRIN2009 grant ``Critical Point Theory and Perturbative Methods for Nonlinear Differential Equations''.}}

\sect Introduction

\sm

One of the main questions in the study of local holomorphic dynamics (see [A1], [A2], [Bra], or [R3] Chapter 1, for general surveys on this topic) is when a germ of biholomorphism~$f$, fixing the origin, is {\it holomorphically linearizable}, i.e., when there exists a local holomorphic change of coordinates such that~$f$ is conjugated to its linear part $\Lambda$. 

A way to solve such a problem is to first look for a formal transformation $\phe$ solving 
$$
f\circ\phe = \phe \circ \Lambda,
$$
i.e., to ask when $f$ is {\it formally linearizable}, and then to check whether $\phe$ is convergent. Moreover, since up to linear changes of the coordinates we can always assume $\Lambda$ to be in Jordan normal form, i.e.,
$$
\Lambda = \pmatrix { \lambda_1 & & \cr
		 			\eps_1   & \lambda_2 & \cr
					 &  \ddots & \ddots & \cr
					 &  & \eps_{n-1} & \lambda_n  },
$$
where the eigenvalues $\lambda_1, \dots, \lambda_n\in\C^*$ are not necessarily distinct, and $\eps_j$ can be non-zero only if $\lambda_{j+1} = \lambda_j$, we can reduce ourselves to study such germs, and to search for $\phe$ {\it tangent to the identity}, that is, with linear part equal to the identity.

\sm The answer to this question depends on the set of eigenvalues of the linear part of~$f$, usually called its {\it spectrum}. In fact, if we denote by~$\Lambda=(\lambda_1, \dots, \lambda_n)\in (\C^*)^n$ the set of the eigenvalues, then it may happen that there exists a multi-index~$Q=(q_1, \dots, q_n)\in \N^n$, with~$|Q|:=\sum_{j=1}^nq_j\ge 2$, such that
$$
\Lambda^Q - \lambda_j:=\lambda_1^{q_1}\cdots\lambda_n^{q_n} - \lambda_j = 0\tag{eqres}
$$
for some~$1\le j\le n$; a relation of this kind is called a {\it (multiplicative) resonance} of~$f$ {\it relative to the $j$-th coordinate}, $Q$ is called a {\it resonant multi-index relative to the $j$-th coordinate}, and we put
$$
{\res}_j(\Lambda):=\{Q\in\N^n\mid |Q|\ge 2, \Lambda^Q = \lambda_j\}.\tag{eqres2}
$$ 
The elements of $\res(\Lambda) := \bigcup_{j=1}^n \res_j(\Lambda)$ are simply called {\it resonant}  multi-indices. A {\it resonant monomial} is a monomial~$z^Q:=z_1^{q_1}\cdots z_n^{q_n}$ in the~$j$-th coordinate with $Q\in\res_j(\Lambda)$.

\sm Resonances are the formal obstruction to linearization. Anyway there are formal, and holomorphic, linearization results also in presence of resonances, see for example [R1] and [R3], and references therein.

\me One generalization of the previous question is to ask when~$h\ge 2$ germs of biholomorphisms~$f_1, \dots, f_h$ of~$\C^n$ at the same fixed point, which we may place at the origin, are {\it simultaneously holomorphically linearizable}, i.e., there exists a local holomorphic change of coordinates conjugating~$f_k$ to its linear part for each~$k=1, \dots, h$.

In dimension $1$, this problem has been thoroughly studied, also for commuting systems of analytic or smooth circle diffeomorphisms, that are indeed deeply related to commuting systems of germs of holomorphic functions, as explained in [P]. The question about the smoothness of a simultaneous linearization of such a system, raised by Arnold, was brilliantly answered by Herman [H], and extended by Yoccoz [Y1] (see also [Y3]). In [M], Moser raised the problem of smooth linearization of commuting circle diffeomorphisms in connection with the holonomy group of certain foliations of codimension 1, and, using the rapidly convergent Nash-Moser iteration scheme, he proved that if the rotation numbers of the diffeomorphisms satisfy a simultaneous Diophantine condition and if the diffeomorphisms are in some $C^\io$-neighborhood of the corresponding rotations (the neighborhood being imposed by the constants appearing in the arithmetic condition, as usual in perturbative KAM theorems) then they are $C^\io$-linearizable, that is, $C^\io$-conjugated to rotations. We refer to [FK] and references therein for a clear exposition of the one-dimensional problem and for the best results, up to now, in such a context. Furthermore, the problem for commuting germs of holomorphic functions in dimension one has been studied by DeLatte [D], and more recently by Biswas [Bi], under Brjuno-type conditions generalizing Moser's simultaneous Diophantine condition.

In dimension $n\ge 2$ much less is know in the formal and holomorphic settings. Gramchev and Yoshino [GY] have proved a simultaneous 
holomorphic linearization result for pairwise commuting germs without simultaneous resonances, with diagonalizable linear parts, and under a simultaneous Diophantine condition (further studied by Yoshino in [Yo]) and a few more technical assumptions. In [DG], DeLatte and Gramchev investigated on holomorphic linearization of germs with linear parts having Jordan blocks, leaving as an open problem the study of simultaneous formal and holomorphic linearization of commuting germs with non-diagonalizable linear parts.
Recently it has been proved in [R2] that $h\ge 2$ germs~$f_1, \dots, f_h$ of biholomorphisms of~$\C^n$, fixing the origin, so that the linear part of~$f_1$ is diagonalizable and~$f_1$ commutes with~$f_k$ for any~$k=2,\dots, h$, under certain arithmetic conditions on the eigenvalues of the linear part of~$f_1$ and some restrictions on their resonances, are simultaneously holomorphically linearizable if and only if there exists a particular complex manifold invariant under~$f_1, \dots, f_h$. 

\eject

Therefore, there are at least three natural questions arising in this setting:

{\narrower
\sm\item{(Q1)} {\it Is it possible to say anything on the shape a (formal) simultaneous linearization can have?}

\sm\item{(Q2)} {\it Are there any conditions on the eigenvalues of the linear parts of $h\ge 2$ germs of simultaneously formally linearizable biholomorphisms ensuring simultaneous holomorphic linearizability?}

\sm\item{(Q3)} {\it Under which conditions on the eigenvalues of the linear parts of $h\ge 2$ pairwise commuting germs of biholomorphisms can one assert the existence of a simultaneous holomorphic linearization of the given germs? In particular, is there a Brjuno-type condition sufficient for convergence?}
\sm

}

Note that the third question is a natural generalization to dimension $n\ge 2$ of the question raised by Moser [M] in the one-dimensional case (see also the introduction of [FK]).

\me  In this paper we shall give complete answers to these three questions without making any assumption on the resonances.
Before stating our answer to the first question, we need the following definition.

\sm\defin{DeAlmostJordanIntro} Let $M_1,\dots,M_h$ be $h\ge 2$ complex $n\times n$ matrices. We say that $M_1,\dots,M_h$ are {\it almost simultaneously Jordanizable}, it there exists a linear change of coordinates $A$ such that $A^{-1} M_1 A,\dots, A^{-1} M_h A$ are {\it almost in simultaneous Jordan normal form}, i.e., for $k=1,\dots,h$ we have
$$
A^{-1} M_k A = \pmatrix { \lambda_{k,1} & & \cr
		 			\eps_{k,1}   & \lambda_{k,2} & \cr
					 &  \ddots & \ddots & \cr
					 &  & \eps_{k,n-1} & \lambda_{k,n}  },~~\eps_{k,j}\ne0\Longrightarrow \lambda_{k,j}=\lambda_{k,j+1}.
\tag{simalmjor}
$$
We say that $M_1,\dots, M_h$ are {\it simultaneously Jordanizable} if there exists a linear change of coordinates $A$ such that we have \rf{simalmjor} with $\eps_{k,j}\in\{0,\eps\}$.

\sm It should be remarked that two commuting matrices are not necessarily
almost simultaneously Jordanizable, and that two almost simultaneously Jordanizable matrices
do not necessarily commute; see section~2 for details. However, the almost simultaneously Jordanizable hypothesis still is less restrictive than the simultaneously 
diagonalizable assumption usual in this context.

\sm The following result gives an answer to (Q1).

\sm\thm{TeLinNotResintro}{Theorem} {\sl Let $f_1, \dots, f_h$ be $h\ge 2$  formally linearizable germs of biholomorphisms of $\C^n$ fixing the origin and with almost simultaneously Jordanizable linear parts. If $f_1, \dots, f_h$ are simultaneously formally linearizable, then they are simultaneously formally linearizable via a linearization $\phe$ such that $\phe_{Q,j} = 0$ for each $Q$ and $j$ so that $Q\in\cap_{k=1}^h \res_j(\Lambda_k)$, and such a linearization is unique.}

\sm We also have a condition ensuring formal simultaneous linearizability.

\sm\thm{TeSimFormLinJordanIntro}{Theorem} {\sl Let $f_1, \dots, f_h$ be $h\ge 2$  formally linearizable germs of biholomorphisms of $\C^n$ fixing the origin and with almost simultaneously Jordanizable linear parts. If $f_1, \dots, f_h$ all commute pairwise, then they are simultaneously formally linearizable.}

\me To state our result on simultaneous holomorphic linearizability we need to introduce the following Brjuno-type condition.

\sm\defin{De1.0bisIntro} Let~$n\ge2$ and let~$\Lambda_1=(\lambda_{1,1}, \dots, \lambda_{1,n}),\dots, \Lambda_h=(\lambda_{h,1}, \dots, \lambda_{h,n})$ be $h\ge2$ $n$-tuples of complex, not necessarily distinct, non-zero numbers. We say that~$\Lambda_1,\dots,\Lambda_h$ {\it satisfy the simultaneous Brjuno condition} if there exists a strictly increasing sequence of integers~$\{p_\nu\}_{\nu_\ge 0}$ with~$p_0=1$ such that
$$
\sum_{\nu\ge 0} {1\over p_\nu} \log{1\over\omega_{\Lambda_1,\ldots,\Lambda_h}(p_{\nu+1})}<+\io,
$$
where for any $m\ge 2$ we set
$$
\omega_{\Lambda_1,\ldots,\Lambda_h}(m)= \min_{2\le|Q|\le m\atop Q\not\in\cap_{k=1}^h\cap_{j=1}^n\res_j(\Lambda_k)}\eps_Q, 
$$
with
$$
\eps_Q =\min_{1\le j \le n}\max_{1\le k\le h}|\Lambda_{k}^Q - \lambda_{k,j}|. 
$$
If $\Lambda_1,\ldots,\Lambda_h$ are the sets of eigenvalues of the linear parts of $f_1,\dots, f_h$, we shall say that $f_1,\dots,f_h$ {\it satisfy the simultaneous Brjuno condition}.

\sm Our holomorphic linearization result answering (Q2) is then the following.

\sm\thm{TeLinSimIntro}{Theorem} {\sl Let $f_1, \dots, f_h$ be $h\ge 2$ simultaneously formally linearizable germs of biholomorphism of $\C^n$ fixing the origin and such that their linear parts $\Lambda_1,\dots, \Lambda_h$ are simultaneously diagonalizable. If $f_1,\dots, f_h$ satisfy the simultaneous Brjuno condition, then $f_1,\dots f_h$ are holomorphically simultaneously linearizable.}

\sm Using Theorem \rf{TeLinSimIntro} we are also able to give a positive answer to the generalization (Q3) of Moser's question.

\sm\thm{TeSimHolLinDiagIntro}{Theorem} {\sl Let $f_1, \dots, f_h$ be $h\ge 2$  formally linearizable germs of biholomorphisms of $\C^n$ fixing the origin, with simultaneously diagonalizable linear parts, and satisfying the simultaneous Brjuno condition. Then $f_1, \dots, f_h$ are simultaneously holomorphically linearizable if and only if they all commute pairwise.}

\me The structure of the paper is as follows. In the next section we shall discuss properties of simultaneously formally linearizable germs, and we shall give a proof of Theorem \rf{TeLinNotResintro}, Theorem~\rf{TeSimFormLinJordanIntro} and other formal results that we can obtain. In section $3$ we shall prove Theorem~\rf{TeLinSimIntro} and Theorem \rf{TeSimHolLinDiagIntro} using majorant series. In the appendix we shall discuss the equivalence between various Brjuno-type series.

\me\no{\bf Acknowledgments.} I would like to thank the members of the Dipartimento di Matematica e Applicazioni of the Universit\`a degli Studi di Milano Bicocca for their warm and nice welcome.

\sect Simultaneously formally linearizable germs

\sm In this section, we shall deal with formal simultaneous linearization. We shall first investigate the properties one can expect from a simultaneous formal linearization, and then we shall provide conditions for simultaneous formal linearizability in presence of resonances. 

\sm\defin{DeUno} Let~$\Lambda\in(\C^*)^n$ and let~$j\in\{1, \dots, n\}$. We say that a multi-index~$Q\in\N^n$, with~$|Q|\ge 2$, gives a {\it resonance relation for~$\Lambda$ relative to the~$j$-th coordinate} if
$$
\Lambda^Q :=\lambda_1^{q_1}\cdots\lambda_n^{q_n} = \lambda_j
$$ 
and we put, as in \rf{eqres2}, ${\res}_j(\Lambda)=\{Q\in\N^n\mid |Q|\ge 2, \lambda^Q = \lambda_j\}$. 
The elements of $\res_j(\Lambda)$ are simply called {\it resonant}  multi-indices with respect to $j$. 

If $\Lambda$ is a complex $n\times n$ invertible matrix in Jordan normal form, then, with a slight abuse of notation, we shall denote by $\res_j(\Lambda)$ the resonant multi-indices of the eigenvalues of $\Lambda$.

\sm Let us start with the following useful result.

\sm\thm{LeJordan}{Lemma} {\sl Let $f$ be a germ of biholomorphism of $\C^n$ fixing the origin, and let $\Lambda$ be an invertible $n\times n$ complex matrix in Jordan normal form, commuting with $f$. Then the linear part of $f$ commutes with $\Lambda$ and $f$ contains only monomials that are resonant with respect to the eigenvalues of $\Lambda$.}

\sm\proof We can write $f$ in coordinates, as $f(z)=Mz +\widehat f(z) =Mz + \sum_{|Q|\ge 2} f_Q z^Q$. If $f$ commutes with $\Lambda$ then, comparing terms of the same degree, it is clear that $M$ has to commute with $\Lambda$.

If $\Lambda$ is diagonal, it is obvious that $f$ commutes with $\Lambda$ if and only if $f$ contains $\Lambda$-resonant terms only.

Let us now assume that $\Lambda$ contains at least a non-trivial Jordan block, that is 
$$
\Lambda = \pmatrix { \lambda_1 & & \cr
		 			\eps_1   & \lambda_2 & \cr
					 &  \ddots & \ddots & \cr
					 &  & \eps_{n-1} & \lambda_n  }, \quad \eps_j\in\{0,\eps\}, ~~\eps_j\ne0\Longrightarrow \lambda_j=\lambda_{j+1},
$$
with at least one non-zero $\eps_j$.

Up to reordering, we may assume $\eps_1\ne 0$.   
For each component $j\in\{1,\dots, n\}$, we have
$$
(\Lambda \widehat f(z) )_j = \lambda_j \widehat f_j(z) + \eps_{j-1} \widehat f_{j-1}(z), 
$$
where we set $\eps_{-1}=0$, and
$$
\eqalign{
(\widehat f(\Lambda z) )_j 
&= \sum_{|Q|\ge 2} f_{Q,j} (\Lambda z)^Q\cr
&= \sum_{|Q|\ge 2} f_{Q,j} \lambda^Q z^Q \prod_{k=2}^n\left(1 + \eps_{k-1}{z_{k-1}\over\lambda_k z_k}\right)^{q_k}\cr
&= \sum_{|Q|\ge 2} f_{Q,j} \lambda^Q z^Q \sum_{{0\le k_2\le q_2\atop\vdots}\atop {0\le k_n\le q_n}} { q_2 \choose k_2}\cdots { q_n  \choose  k_n} {\eps_1^{k_2}\cdots\eps_{n-1}^{k_n}\over \lambda_2^{k_2}\cdots \lambda_n^{k_n} } z_1^{k_2} z_2^{k_3-k_2}\cdots z_{n-1}^{k_n-k_{n-1}}z_n^{-k_n}\;,
}
$$
where we are using the convention $0^0 = 1$.  
Note that $z^P = z^Q\cdot z_1^{k_2} z_2^{k_3-k_2}\cdots z_{n-1}^{k_n-k_{n-1}}z_n^{-k_n}$ is a monomial with the same degree as $z^Q$, i.e., $|P| = |Q|$,  but subsequent in the lexicographic order, i.e., $P>Q$. Moreover, if $Q$ is a resonant multi-index relative to $j$, i.e., $\lambda^Q=\lambda_j$, then also $P$ is. In fact, for each Jordan block of order $\ell$, we have $\lambda_{i_\ell}=\cdots=\lambda_{i_\ell + \ell}$, so 
$$
\lambda_{i_\ell}^{q_{i_\ell}}\cdots\lambda_{i_\ell + \ell}^{q_{i_\ell + \ell}}= \lambda_{i_\ell}^{q_{i_\ell}+k_{i_\ell+1}}\lambda_{i_\ell+1}^{q_{i_\ell+1}+k_{i_\ell+2}-k_{i_\ell+1}}\cdots\lambda_{i_\ell + \ell}^{q_{i_\ell + \ell}-k_{i_\ell +\ell}}.
$$
Now we prove that for each $j\in\{1,\dots, n\}$ the $j$-th component of $\widehat f$ contains only $\Lambda$-resonant monomials. 
For the first component we have 
$$
\lambda_1\widehat f_1(z) = \sum_{|Q|\ge 2} f_{Q,1} \lambda^Q z^Q \sum_{{0\le k_2\le q_2\atop\vdots}\atop {0\le k_n\le q_n}} { q_2 \choose k_2}\cdots { q_n  \choose  k_n} {\eps_1^{k_2}\cdots\eps_{n-1}^{k_n}\over \lambda_2^{k_2}\cdots \lambda_n^{k_n} } z_1^{k_2} z_2^{k_3-k_2}\cdots z_{n-1}^{k_n-k_{n-1}}z_n^{-k_n}. 
\tag{eqf1}
$$ 
Let $\widetilde Q$ be the first, with respect to the lexicographic order, non resonant multi-index so that $f_{\widetilde Q, 1}\ne 0$ and let us compare the coefficients of $z^{\widetilde Q}$ in both sides of \rf{eqf1}. In the left-hand side we just have $\lambda_1 f_{\widetilde Q, 1}$; in the right-hand side we only have $\lambda^{\widetilde Q} f_{\widetilde Q, 1}$, because other contributes could come only by previous multi-indices, but, as observed above, they all give resonances because we are assuming $\widetilde Q$ to be the first non resonant multi-index. Hence we have
$$
(\lambda^{\widetilde Q} -\lambda_1)f_{\widetilde Q, 1} = 0, 
$$
yielding, since $\lambda^{\widetilde Q} \ne\lambda_1$, $f_{\widetilde Q, 1} = 0$, and contradicting the hypothesis.
Now we turn to the second component, and since we are assuming $\eps_1\ne0$, we have $\lambda_2=\lambda_1$, so we have
$$
\lambda_1\widehat f_2(z) +\eps \widehat f_1(z) \!= \!\!\!\sum_{|Q|\ge 2} \!\!\!f_{Q,2} \lambda^Q z^Q \!\!\!\!\!\sum_{{0\le k_2\le q_2\atop\vdots}\atop {0\le k_n\le q_n}} \!\!\!\!\!{ q_2 \choose k_2}\!\cdots \!{ q_n  \choose  k_n} {\eps_1^{k_2}\cdots\eps_{n-1}^{k_n}\over \lambda_2^{k_2}\cdots \lambda_n^{k_n} } z_1^{k_2} z_2^{k_3-k_2}\cdots z_{n-1}^{k_n-k_{n-1}}z_n^{-k_n}. 
\tag{eqf2}
$$
Let $\widetilde Q$ be the first, with respect to the lexicographic order, non resonant multi-index so that $f_{\widetilde Q, 2}\ne 0$ and let us compare the coefficients of $z^{\widetilde Q}$ in both sides of \rf{eqf2}. In the left-hand side we just have $\lambda_1 f_{\widetilde Q, 2}$ because we proved above that $\widehat f_1$ contains only resonant monomials; in the right-hand side we again have only $\lambda^{\widetilde Q} f_{\widetilde Q, 2}$, because other contributes could come only by previous multi-indices, but, as observed above, they all give resonances since $\widetilde Q$ is the first non resonant multi-index. Then we repeat the same argument used above and we prove that also $\widehat f_2$ contains only resonant monomials. Now we can use the same arguments for the remaining components, and we get the assertion. \qed

\sm\thm{ReLemma0}{Remark} Notice that in the previous result we did not make any hypotheses on the diagonalizability or not of the linear part $M$ of the germ $f$, because we just wanted to understand what information we can deduce on $f$ assuming its commutation with a matrix in Jordan normal form.

\sm\thm{ReLemma}{Remark} Note that Lemma \rf{LeJordan} does not hold if $\Lambda$ is just triangular. For example, if we take
$$
\Lambda=\pmatrix{\lambda_2\lambda_3&0&0\cr 0&\lambda_2&0\cr \lambda_3(1-\lambda_2)&
\lambda_3-\lambda_2&\lambda_3\cr}
$$
with $\lambda_2,\lambda_3\in\C\setminus\{0,1\}$ and $\lambda_2\ne \lambda_3$, and $A$ is any complex $3\times3$ matrix commuting with $\Lambda$ (for example $\Lambda$ itself), then $\Lambda$ commutes with the germ
$$
f(z)=Az+\left(z_2(z_1+z_2+z_3), 0, -z_2(z_1+z_2+z_3)\right),
$$ 
and $f$ clearly contains monomials non resonant with respect to the eigenvalues of $\Lambda$.

\bi In the following we shall need the notion of almost simultaneously Jordanizable defined in the introduction in Definition \rf{DeAlmostJordanIntro}, that we recall here.

\sm\defin{DeAlmostJordan} Let $M_1,\dots,M_h$ be $h\ge 2$ complex $n\times n$ matrices. We say that $M_1,\dots,M_h$ are {\it almost simultaneously Jordanizable}, it there exists a linear change of coordinates $A$ such that $A^{-1} M_1 A,\dots, A^{-1} M_h A$ are {\it almost in simultaneous Jordan normal form}, i.e., for $k=1,\dots,h$ we have
$$
A^{-1} M_k A = \pmatrix { \lambda_{k,1} & & \cr
		 			\eps_{k,1}   & \lambda_{k,2} & \cr
					 &  \ddots & \ddots & \cr
					 &  & \eps_{k,n-1} & \lambda_{k,n}  },~~\eps_{k,j}\ne0\Longrightarrow \lambda_{k,j}=\lambda_{k,j+1}.
\tag{simalmjor}
$$
We say that $M_1,\dots, M_h$ are {\it simultaneously Jordanizable} if there exists a linear change of coordinates $A$ such that we have \rf{simalmjor} with $\eps_{k,j}\in\{0,\eps\}$.

\sm\thm{ReJordanizable}{Remark} Note that the problem of deciding when two $n\times n$ complex matrices are almost simultaneously Jordanizable is not as easy as when the two matrices are diagonalizable. Indeed, whereas $h\ge 2$ diagonalizable matrices are simultaneously diagonalizable if and only if they commute pairwise, and if $h\ge 2$ matrices commute pairwise then they are simultaneously triangularizable (but the converse is clearly false), if 
two matrices commute then this does not imply that they admit an almost simultaneous Jordan normal form, and it is not true in general that two matrices almost in simultaneous Jordan normal form commute. For example the following two matrices
$$
\Lambda= \pmatrix{\lambda & 0 & 0\cr	
				  \eps & \lambda & 0\cr
				  0& 0& \lambda}
\quad	
M= \pmatrix{\mu & 0 & 0\cr	
			\delta & \mu & 0\cr
			\beta & 0& \mu} \quad \lambda, \eps, \mu,\delta,\beta\in\C^*
$$
commute, but in general it is not possible to almost simultaneously Jordanize them. In fact all the matrices $A$ such that $M$ is almost in simultaneous Jordan normal form with $\Lambda$ have to be invertible solutions of the following equation
$$
AM=\pmatrix{\mu&0&0\cr \zeta&\mu&0\cr 0&0&\mu\cr} A
$$
(and $\zeta$ has to be non-zero because $M\ne\mu I_3$); hence $A$ is of the form
$$
A=\pmatrix{{\beta\over\zeta}f+{\delta\over\zeta}e&0&0\cr d&e&f\cr g&h&-{\delta\over\beta}h\cr}\;,
$$
which is invertible if and only if $\beta f+\delta e\ne 0$ and $h\ne 0$. But
$$
A\Lambda=\pmatrix{\lambda&0&0\cr \xi&\lambda&0\cr 0&0&\lambda\cr} A
$$
yields $h=0$, implying that $\Lambda$ and $M$ are not almost simultaneously Jordanizable (and so they are also not simultaneously Jordanizable).

On the other side, the following two matrices
$$
\widetilde\Lambda= \pmatrix{\lambda & 0 & 0\cr	
				  \eps & \lambda & 0\cr
				  0& \eps& \lambda}
\quad	
\widetilde M= \pmatrix{\mu & 0 & 0\cr	
			\delta & \mu &0 \cr
			0 & 0& \eta} \quad \lambda, \eps, \mu,\delta,\eta\in\C^*
$$
are almost in simultaneous Jordan normal form but they do not commute. 

\me In studying the convergence of a formal linearization of a germ of biholomorphism, and hence also in the case of simultaneous formal linearizations, it is useful to be able to use formal linearizations of the special kind we are now going to introduce. 

\sm\defin{DeLinearNonRes} Let $f$ be a germ of bihomolomorphism of $\C^n$ fixing the origin and with linear part $\Lambda$ in Jordan normal form. A tangent to the identity (formal) linearization $\phe$ of $f$ is said {\it non resonant} if for each resonant multi-index relative to the $j$-th coordinate, $Q\in\res_j(\Lambda)$, the coefficient $\phe_{Q,j}$ of $z^Q$ in the power series expansion of the $j$-th coordinate of $\phe$ vanishes, i.e., $\phe_{Q,j} = 0$.

\sm\defin{DeSimLinearNonRes} Let $f_1, \dots, f_h$ be $h\ge 2$  formally linearizable germs of biholomorphisms of $\C^n$ fixing the origin and with linear parts $\Lambda_1, \dots, \Lambda_h$ almost in simultaneous Jordan normal form. A tangent to the identity (formal) simultaneous linearization $\phe$ of $f_1, \dots, f_h$ is said {\it non resonant} if for each simultaneous resonant multi-index relative to the $j$-th coordinate, $Q\in\cap_{k=1}^h\res_j(\Lambda_k)$, the coefficient $\phe_{Q,j}$ of $z^Q$ in the power series expansion of the $j$-th coordinate of $\phe$ vanishes, i.e., $\phe_{Q,j} = 0$.

\sm Let us now investigate the shape a formal simultaneous linearization can have. It can be proven, see [R\"u] and [R4], that a formally linearizable germ of biholomorphism is formally linearizable via a non resonant formal linearization, and such a linearization is unique. The same is true also for simultaneously formally linearizable germs, but with a slightly different proof, as shown in the next result, stated in the introduction as Theorem \rf{TeLinNotResintro}.

\sm\thm{TeLinNotRes}{Theorem} {\sl Let $f_1, \dots, f_h$ be $h\ge 2$  formally linearizable germs of biholomorphisms of $\C^n$ fixing the origin and with almost simultaneously Jordanizable linear parts. If $f_1, \dots, f_h$ are simultaneously formally linearizable, then they are simultaneously formally linearizable via a non resonant linearization $\phe$, and such a linearization is unique.}

\sm\proof We may assume, up to linear changes of the coordinates, that the linear parts $\Lambda_1, \dots, \Lambda_h$ of $f_1, \dots, f_h$ are almost in simultaneous Jordan normal form, i.e.,
$$
\Lambda_k = \pmatrix { \lambda_{k,1} & & \cr
		 			\eps_{k,1}   & \lambda_{k,2} & \cr
					 &  \ddots & \ddots & \cr
					 &  & \eps_{k,n-1} & \lambda_{k,n}  },~~\eps_{k,j}\ne0\Longrightarrow \lambda_{k,j}=\lambda_{k,j+1},
$$
for $k=1,\dots, h$. 

We know that there exists a formal change of coordinates $\phe$ tangent to the identity and such that $\phe^{-1}\circ f_k\circ \phe = \Lambda_k$ for all $k=1,\dots, h$. If $\phe_{Q,j} = 0$ for each $Q$ and $j$ so that $Q\in\cap_{k=1}^h \res_j(\Lambda_k)$, then we are done and we only have to show that such a linearization is unique.

If there is at least one multi-index $Q\in\cap_{k=1}^h \res_j(\Lambda_k)$ with $j\in\{1,\dots, n\}$ and such that $\phe_{Q,j}\ne 0$, then we can construct another formal simultaneous linearization $\psi$ which is non-resonant. Since we can write $f_k=\phe\circ \Lambda_k\circ\phe^{-1}$ for each $k=1,\dots, h$, $\psi$ has to satisfy
$$
\phe^{-1}\circ \psi\circ \Lambda_k = \Lambda_k\circ\phe^{-1}\circ \psi,
$$
i.e., we need to construct $\psi$ not containing monomials simultaneously resonant for $\Lambda_1, \dots, \Lambda_h$ and such that (by Lemma \rf{LeJordan}) $\phe^{-1}\circ \psi$ contains only monomials simultaneously resonant for $\Lambda_1,\dots, \Lambda_h$. 
Writing 
$$
\phe^{-1}_j(z) = z_j \left(1 + \sum_{Q\in N_j \atop \Lambda_k^Q = 1, k=1,\dots, h} \widetilde\phe_{Q,j}z^Q + \sum_{Q\in N_j \atop \Lambda_k^Q \ne 1,~\rm{for~some}~k=1,\dots, h} \widetilde\phe_{Q,j}z^Q\right), \quad j=1,\dots, n,
$$
where
$$
N_j:=\{Q\in \Z^n \mid |Q|\ge 1, q_j\ge -1, q_h\ge 0~\hbox{for all}~h\ne j\},
$$
we want to find $\psi$ of the form
$$
\psi_j(z) = z_j\left(1 + \sum_{Q\in N_j \atop \Lambda_k^Q \ne 1, ~\rm{for~some}~k=1,\dots, h} \psi_{Q,j}z^Q\right),\quad j=1,\dots, n,
$$
such that 
$$
(\phe^{-1}\circ\psi)_j(z) = z_j\left(1 + \sum_{Q\in N_j \atop \Lambda_k^Q = 1, k=1,\dots, h} g_{Q,j}z^Q\right), \quad j=1,\dots, n.
$$
Using the same argument of the proof of Poincar\'e-Dulac Theorem (see [Ar] pp.~192--193 or [R3] pp. 40--41), in the $j$-th coordinate of $\phe^{-1}\circ \psi$, the coefficient of $z^Q$ with $Q$ so that $\Lambda_k^Q\ne 1$ for at least one $k$ is equal to
$$
\psi_{Q,j} + \hbox{Polynomial}(\phe^{-1}, \hbox{previous}~\psi_{P, l}),
$$
where the polynomial in the previous formula is ``universal'' in the sense that it depends only on the fact that we are composing two power series and it does not depend on the coefficients of $\phe^{-1}$ and on the previous $\psi_{P,l}$, that are in fact just arguments of this universal polynomial.

\no Hence it suffices to put 
$$
\psi_{Q,j} =- \hbox{Polynomial}(\phe^{-1}, \hbox{previous}~\psi_{P, l}).
$$
Note that for the first non-resonant multi-indices, due to degree considerations, we just have to put $\psi_{Q,j} = -\widetilde \phe_{Q,j}$.

We proved that there exists a formal non-resonant tangent to the identity simultaneous linearization $\phe$ of the given germs, containing only monomials that are not simultaneous resonant for the eigenvalues of $\Lambda_1,\dots,\Lambda_h$. Let us assume by contradiction that there exists another such a linearization $\psi$ and $\psi\not \equiv \phe$. Writing $\phe(z) = z+ \sum_{|Q|\ge 2}\phe_Q z^Q$ and $\psi(z) = z+ \sum_{|Q|\ge 2}\psi_Q z^Q$, let $\widetilde Q$ be the first multi-index, with respect to the lexicographic order, so that $\phe_{\widetilde Q}\ne\psi_{\widetilde Q}$ and let $\ell\in\{1,\dots, n\}$ be the minimal index such that $\phe_{\widetilde Q,\ell}\ne\psi_{\widetilde Q,\ell}$. Since for each $Q\in\cap_{k=1}^h \res_\ell (\Lambda_k)$ we know that $\phe_{Q,\ell}=\psi_{Q,\ell}=0$, there is at least one germ $f_k$ such that $\Lambda_k^{\widetilde Q}\ne\lambda_{k,\ell}$.
We know, again by the proof of Poincar\'e-Dulac Theorem, that
$$
\left(\Lambda_k^{\widetilde Q} -\lambda_{k,\ell}\right) \phe_{\widetilde Q,\ell} = \hbox{Polynomial}(f_k, \hbox{previous}~\phe_{P, j}),
$$
and, similarly,
$$
\left(\Lambda_k^{\widetilde Q} -\lambda_{k,\ell}\right) \psi_{\widetilde Q,\ell} = \hbox{Polynomial}(f_k, \hbox{previous}~\psi_{P, j}),
$$
where, again, the polynomial in the previous formulas is ``universal'', because again it depends only on the fact that we are composing power series and it does not depend on the coefficients of $f_k$ and on the previous $\phe_{P,j}$ or $\psi_{P,j}$, that are in fact just arguments of this universal polynomial. Hence, since we are assuming that $\phe_{P, j}= \psi_{P,j}$ for all the multi-indices $P<\widetilde Q$ and, for $P=\widetilde Q$, for all $j< \ell$, we have that 
$$
\hbox{Polynomial}(f_k, \hbox{previous}~\phe_{P, j})=\hbox{Polynomial}(f_k, \hbox{previous}~\psi_{P, j})
$$
implying that $\phe_{\widetilde Q,\ell}=\psi_{\widetilde Q,\ell}$, and contradicting the hypothesis. \qed

\sm\thm{ReEcalle}{Remark} The universal polynomials we dealt with in the last proof can be interpreted and computed using the mould formalism introduced by \'Ecalle (see [\'E]), and in this sense the latter proof is a mouldian proof. 

\sm As announced in the introduction as Theorem \rf{TeSimFormLinJordanIntro}, we shall now give a condition ensuring formal simultaneous linearizability. 

\sm\thm{TeSimFormLinJordan}{Theorem} {\sl Let $f_1, \dots, f_h$ be $h\ge 2$  formally linearizable germs of biholomorphisms of $\C^n$ fixing the origin and with almost simultaneously Jordanizable linear parts. If $f_1, \dots, f_h$ all commute pairwise, i.e., $f_p\circ f_q = f_q\circ f_p$ for any $p$ and $q$ in $\{1,\dots, h\}$, then they are simultaneously formally linearizable.}

\sm\proof We may assume without loss of generality that the linear parts $\Lambda_1,\dots, \Lambda_h$ of the germs are all almost in simultaneous Jordan normal form. Since $f_1$ is formally linearizable, it is possibly to linearize it with a  non-resonant formal linearization $\phe_1$ (see [R\"u] and [R4]). Then, thanks to the commutation hypothesis, by Lemma \rf{LeJordan}, $\widetilde f_2 = \phe_1^{-1}\circ f_2\circ \phe_1, \dots,\widetilde f_h = \phe_1^{-1}\circ f_h\circ \phe_1 $ contain only $\Lambda_1$-resonant terms. Now we claim that it is possible to find a formal change of coordinates $\phe_2$ fixing the origin, tangent to the identity, containing only $\Lambda_1$-resonant terms that are not $\Lambda_2$-resonant, and conjugating $\widetilde f_2$ to a germ $g_2$ with same linear part and in Poincar\'e-Dulac normal form. In fact, we have to solve
$$
\left(\Lambda_2 + \widetilde f_2^{\res(\Lambda_1)}\right) \circ \left(I + \phe_2^{\res(\Lambda_1)}\right) =  
\left(I + \phe_2^{\res(\Lambda_1)}\right) \circ \left(\Lambda_2 + g_2^{\res(\Lambda_1)\cap\res(\Lambda_2)}\right), 
$$
where $\widetilde f^{\res(\Lambda_1)}$ contains only monomial resonant with respect to the eigenvalues of $\Lambda_1$, and so on, that is
$$
\Lambda_2\phe_2^{\res(\Lambda_1)} +  \widetilde f_2^{\res(\Lambda_1)} \circ \left(I + \phe_2^{\res(\Lambda_1)}\right) =  
g_2^{\res(\Lambda_1)\cap\res(\Lambda_2)} + \phe_2^{\res(\Lambda_1)} \circ \left(\Lambda_2 + g_2^{\res(\Lambda_1)\cap\res(\Lambda_2)}\right), 
$$
which is solvable by the usual Poincar\'e-Dulac procedure and setting $\phe_{2, Q, j}^{\res(\Lambda_1)}=0$ whenever $\Lambda_2^Q=\lambda_{2,j}$. But since $f_2$ is formally linearizable, the linear form is its unique Poincar\'e-Dulac normal form (see [R4] Theorem 2.3), so $g_2\equiv \Lambda_2$. 
Hence, since the given germs commute pairwise, $\phe_2^{-1}\circ\phe_1^{-1}\circ f_3\circ\phe_1\circ\phe_2, \dots, \phe_2^{-1}\circ\phe_1^{-1}\circ f_h\circ\phe_1\circ\phe_2$ contain only monomials that are simultaneously $\Lambda_1$ and $\Lambda_2$ resonant, and we can iterate the procedure finding a formal linearization $\phe_3$ of $\phe_2^{-1}\circ\phe_1^{-1}\circ f_3\circ\phe_1\circ\phe_2$ containing only monomials that are $\Lambda_1$ and $\Lambda_2$ resonant but not $\Lambda_3$-resonant. We can then iteratively perform the same procedure getting $\phe_4,\dots,\phe_h$ formal transformations such that $\phe_1\circ\cdots\circ\phe_h$ is a simultaneous formal linearization of $f_1,\dots,f_h$. \qed

\sm\thm{ReSimLinForm}{Remark} The simultaneous formal linearization obtained in the last proof is non-resonant.
Moreover, we can perform the same procedure with a different permutation of the indices, i.e., starting with $f_{\sigma(1)}$ and then continuing with $f_{\sigma(2)}$ and so on, where $\sigma$ is any permutation of $\{1,\dots,n\}$, and, by Theorem \rf{TeLinNotRes}, we always get the same linearization.

\sm\thm{ReSimLin}{Remark} The hypothesis on the pairwise commutation is indeed necessary. In fact, if $\Lambda_1$ and $\Lambda_2$ are two commuting matrices almost in simultaneous Jordan normal form such that $\res(\Lambda_1) \ne\void$ and $\res(\Lambda_2) \ne \void$, but $\res(\Lambda_1)\cap\res(\Lambda_2)=\void$, the unique formal transformation tangent to the identity and commuting with both $\Lambda_1$ and $\Lambda_2$ is the identity, so any non-linear germ $f_3$ with linear part in Jordan normal form and commuting with $\Lambda_1$ (that is, containing only $\Lambda_1$-resonant terms) but not with $\Lambda_2$ cannot be simultaneously linearizable with $\Lambda_1$ and $\Lambda_2$. 

\me Note that, in the proofs of the previous results, we needed to assume the almost simultaneous Jordanizability of the linear parts of the given germs, and indeed we cannot perform those proofs just assuming that those linear parts are just simultaneously triangularizable, because, as already remarked, Lemma \rf{LeJordan} does not hold for general triangular matrices. Anyway, this is more than what was usually known in the previous literature, where in general linearization results are proved only for germs with diagonalizable linear part.
Moreover, in the particular case of diagonalizable linear parts, recalling that $h\ge 2$ diagonalizable commuting complex $n\times n$ matrices are simultaneously diagonalizable, we have the following equivalence. 

\sm\thm{TeSimFormLinDiag}{Proposition} {\sl Let $f_1, \dots, f_h$ be $h\ge 2$  formally linearizable germs of biholomorphisms of $\C^n$ fixing the origin and with simultaneously diagonalizable linear parts. Then $f_1, \dots, f_h$ are simultaneously formally linearizable if and only if they all commute pairwise.}

\sm\proof If $f_1, \dots, f_h$ are simultaneously formally linearizable, since their linear parts are simultaneously diagonalizable, then they all commute and we are done.

The converse follows from Theorem \rf{TeSimFormLinJordan}. \qed

\sm\thm{ReLinNoCommut}{Remark} Note that having simultaneously diagonalizable linear parts is equivalent to having diagonalizable pairwise commuting linear parts, but being simultaneously formally linearizable does not imply that the linear parts are pairwise commuting, even when the linear parts are diagonalizable.

\me It is possible to find cases of formal simultaneous linearization even without assuming that the linear parts are almost simultaneously Jordanizable or that the germs commute pairwise, 
as shown in the following results. However, in those cases one has to assume other conditions, for example on the nature of resonances. We refer to [R2] and [R3] for the definitions of {\it only level $s$ resonances} and {\it simultaneous osculating manifold}. 

\sm\thm{PropOsculanteFormale}{Proposition} {\sl Let $f_1,\dots,f_h$ be $h\ge 2$ germs of biholomorphism of $\C^n$ fixing the origin. Assume that the spectrum of the linear part of $f_1$ has only level $s$ resonances and that $f_1$ commutes with $f_k$ for $k=2,\dots, h$. Then $f_1,\dots, f_h$ are simultaneously formally linearizable if and only if there exists a germ of formal complex manifold~$M$ at~$O$ of codimension~$s$, invariant under~$f_h$ for each~$h=1, \dots, m$, which is a simultaneous osculating manifold for~$f_1, \dots, f_m$ and such that~$f_1|_M, \dots, f_m|_M$ are simultaneously formally linearizable.}

\sm\proof It is clear from the proof of Theorem $2.5$ of [R2] that it is possible to perform it in the formal category, where to formally linearize $f_1$, under these hypotheses, it suffices to have the only level $s$ hypothesis. \qed

\sm\thm{LeSimFormLin}{Proposition} {\sl Let $f_1, \dots, f_h$ be $h\ge 2$ germs of biholomorphisms of $\C^n$ fixing the origin and such that the eigenvalues of the linear part of $f_1$ have no resonances. If $f_1$ commutes with $f_k$ for $k=2,\dots, h$, then $f_1,\dots, f_h$ are simultaneously formally linearizable via a unique formal linearization.}

\sm\proof Up to linear conjugacy, we may assume that the linear part $\Lambda_1$ of $f_1$ is in Jordan normal form. Since the eigenvalues of $\Lambda_1$ are non-resonant, $f_1$ is formally linearizable via a unique linearization $\phe$ with no constant term and with identity linear part. Then, by Lemma \rf{LeJordan}, thanks to the commutation hypothesis, $\phe^{-1}\circ f_k\circ \phe$ contains only monomials that are resonant with respect to the eigenvalues of $\Lambda_1$, and we get the assertion. \qed

\sm In particular, if the eigenvalues of the linear part of all $f_k$ have no resonances, the unique formal linearization $\phe_k$ of $f_k$ is indeed the unique simultaneous formal linearization of the germs $f_1,\dots, f_h$. 

\sm Another example of the results we can obtain using this kind of arguments is the following. 

\sm\thm{CoUno}{Corollary} {\sl Let $f_1, \dots, f_h$ be $h\ge2$ germs of biholomorphisms of $\C^n$ fixing the origin and such that the linear parts $\Lambda_1$ and $\Lambda_2$ of $f_1$ and $f_2$ are almost simultaneously Jordanizable and have no common resonances (i.e., $\res_j(\Lambda_1) \cap \res_j(\Lambda_2) = \void$ for $j=1,\dots, n$). If $f_1$ and $f_2$ commute with $f_k$ for $k=1,\dots, h$, then $f_1,\dots, f_h$ are simultaneously formally linearizable.}

\sect Convergence under the simultaneous Brjuno condition

\sm In this section we shall prove a holomorphic simultaneous linearization result in presence of resonances. Note that, even if in the previous section we proved formal simultaneous linearization results just assuming almost simultaneous Jordanizability of the linear parts, to prove convergence in ``generic cases'' we need to assume simultaneous diagonalizability of the linear parts because of Yoccoz's counter example (see [Y2] pp. 83--85) to holomorphic
linearization in case of linear part in (non diagonal) Jordan form.

\sm We shall now recall the arithmetical Brjuno-type condition, stated in the introduction in Definition \rf{De1.0bisIntro}, that we shall need to prove our holomorphic linearization result.

\sm\defin{De1.0} Let~$n\ge2$ and let~$\Lambda_1=(\lambda_{1,1}, \dots, \lambda_{1,n}),\dots, \Lambda_h=(\lambda_{h,1}, \dots, \lambda_{h,n})$ be $h\ge2$ $n$-tuples of complex, not necessarily distinct, non-zero numbers. For any $m\ge 2$ we set
$$
\omega_{\Lambda_1,\ldots,\Lambda_h}(m)= \min_{2\le|Q|\le m\atop Q\not\in\cap_{k=1}^h\cap_{j=1}^n\res_j(\Lambda_k)}\eps_Q, 
$$
with
$$
\eps_Q =\min_{1\le j \le n}\max_{1\le k\le h}|\Lambda_{k}^Q - \lambda_{k,j}|. 
$$

Note that, since $Q\not\in\cap_{k=1}^h\cap_{j=1}^n\res_j(\Lambda_k)$, we have $\eps_Q$ always non-zero, and hence also $\omega(m)$ is always non-zero. 
If $\Lambda_1,\ldots,\Lambda_h$ are the sets of eigenvalues of the linear parts of~$f_1, \dots, f_h$, we shall write~$\omega_{f_1,\dots,f_h}(m)$ for~$\omega_{\Lambda_1,\ldots,\Lambda_h}(m)$.

\sm\defin{De1.0bis} Let~$n\ge2$ and let~$\Lambda_1=(\lambda_{1,1}, \dots, \lambda_{1,n}),\dots, \Lambda_h=(\lambda_{h,1}, \dots, \lambda_{h,n})$ be $h\ge2$ $n$-tuples of complex, not necessarily distinct, non-zero numbers. We say that~$\Lambda_1,\dots,\Lambda_h$ {\it satisfy the simultaneous Brjuno condition} if there exists a strictly increasing sequence of integers~$\{p_\nu\}_{\nu_\ge 0}$ with~$p_0=1$ such that
$$
\sum_{\nu\ge 0} {1\over p_\nu} \log{1\over\omega_{\Lambda_1,\ldots,\Lambda_h}(p_{\nu+1})}<+\io.
\tag{simBrjuno}
$$
If $\Lambda_1,\ldots,\Lambda_h$ are the sets of eigenvalues of the linear parts of~$f_1, \dots, f_h$, we shall say that $f_1,\dots,f_h$ {\it satisfy the simultaneous Brjuno condition}.

\sm Recall that \rf{simBrjuno} is equivalent to
$$
\sum_{\nu\ge 0} {1\over 2^\nu} \log{1\over\omega_{\Lambda_1,\ldots,\Lambda_h}(2^{\nu+1})}<+\io,
$$
and a proof of this equivalence can be found in [Brj] and [R3]. 

\sm We shall give in the appendix other conditions equivalent to the simultaneous Brjuno condition. 

\sm\thm{ReAppendice}{Remark} Note that if $\Lambda_1,\dots,\Lambda_h$ do not satisfy the simultaneous Brjuno condition, then each of them does not satisfy the reduced Brjuno condition, i.e.,
$$
\sum_{\nu\ge 0} {1\over 2^\nu} \log{1\over\omega_{\Lambda_k}(2^{\nu+1})}= +\io,
$$
for $k=1,\dots, h$, where 
$$
\omega_{\Lambda_k}(m):= \min_{{2\le|Q|\le m\atop 1\le j \le n}\atop Q\not\in\res_j(\Lambda_k) }|\Lambda_{k}^Q - \lambda_{k,j}|.
$$
In particular, if $\Lambda_1,\dots,\Lambda_h$ are simultaneously Cremer, i.e.,
$$
\limsup_{m\to+\infty}{1\over
m}\log{1\over\omega_{\Lambda_1,\ldots,\Lambda_h}(m)}=+\infty,
$$
and hence they do not satisfy the simultaneous Brjuno condition, then at least one of them has to be Cremer, i.e., 
$$
\limsup_{m\to+\infty}{1\over
m}\log{1\over\omega_{\Lambda_k}(m)}=+\infty,
\tag{eqCremer}
$$
and the other ones do not satisfy the reduced Brjuno condition.

\sm Furthermore it is possible to find $\Lambda_1,\dots,\Lambda_h$ satisfying the simultaneous Brjuno condition, with $\Lambda_k$ not satisfying the reduced Brjuno condition for any $k=1,\dots, h$, as shown in the next result, descending from Theorem 5.1 of [Yo] and Theorem 2.1 of [M].

\sm\thm{PrYoshinoMoser}{Proposition} {\sl Let $h>n\ge2$. Then there exists a set of $h$-tuples of linearly independent vectors $\Lambda_1,\dots, \Lambda_h\in(\C^*)^n$, with the density of continuum, satisfying the simultaneous Brjuno condition, whereas for any $p_1,\dots,p_h\in\Z\setminus\{0\}$ the vector $\sum_{k=1}^h p_k\Lambda_k$ does not satisfy the reduced Brjuno condition.}

\sm\proof The result descends easily from the proof of Theorem 5.1 of [Yo]. In fact, since
$$
\min_{1\le j\le n\atop Q\not \in\cap_{k=1}^h\res_j(\Lambda_k)}\sum_{k=1}^h|\Lambda_k^Q-\lambda_{k,j}| = \sum_{k=1}^h\min_{1\le j\le n\atop Q\not \in\cap_{k=1}^h\res_j(\Lambda_k)}\!\!\!\!\!\!\!\!\!\!|\Lambda_k^Q-\lambda_{k,j}|\le h \cdot\!\!\!\!\!\!\!\!\!\min_{1\le j\le n\atop Q\not \in\cap_{k=1}^h\res_j(\Lambda_k)}\!\!\max_{1\le k\le h}|\Lambda_k^Q-\lambda_{k,j}|,
$$
it is clear that the simultaneous Diophantine condition introduced by Yoshino in Section 5 of [Yo] implies our simultaneous Brjuno condition. 

Hence we can use exactly the same construction of Yoshino in the proof of Theorem 5.1 of [Yo] observing that, since the set of Cremer points, i.e., of the numbers $\theta\in\R$ satisfying
$$
\limsup_{m\to+\infty}{1\over
m}\log{1\over\min_{2\le k\le m}|e^{2\pi i k\theta}-1|}=+\infty,
$$
is residual, i.e., countable intersection of open dense sets, we can perform the same argument substituting, at the end of the proof of Theorem $5.1$ of [Yo], Liouville with Cremer, and we are done. \qed

\sm\thm{RmkSimBrjuno}{Remark} Other possible definitions for a simultaneous Brjuno condition could had been given using the following functions:
$$
\overline\omega_{\Lambda_1,\ldots,\Lambda_h} (m):=\min_{2\le|Q|\le m\atop Q\not\in\cap_{k=1}^h\cap_{j=1}^n\res_j(\Lambda_k)}\max_{1\le k\le h}\min_{1\le j \le n}|\Lambda_{k}^Q - \lambda_{k,j}|, 
$$
or 
$$
\widetilde\omega_{\Lambda_1,\ldots,\Lambda_h} (m) := \max_{1\le k\le h} \min_{2\le|Q|\le m\atop Q\not\in\cap_{k=1}^h\cap_{j=1}^n\res_j(\Lambda_k)}\min_{1\le j \le n}|\Lambda_{k}^Q - \lambda_{k,j}|,
$$
and asking for 
$$
\sum_{\nu\ge 0} {1\over 2^\nu} \log\overline\omega_{\Lambda_1,\ldots,\Lambda_h} (2^{\nu+1})<+\io,
$$
or
$$
\sum_{\nu\ge 0} {1\over 2^\nu} \log\widetilde\omega_{\Lambda_1,\ldots,\Lambda_h} (2^{\nu+1})<+\io.
$$
However, this definitions are clearly more restrictive than the one we chose, since we have:
$$
\widetilde\omega_{\Lambda_1,\ldots,\Lambda_h} (m) \le \overline\omega_{\Lambda_1,\ldots,\Lambda_h} (m)\le \omega_{\Lambda_1,\ldots,\Lambda_h} (m),
$$
and hence
$$
\sum_{\nu\ge 0} {1\over 2^\nu} \log{1\over\omega_{\Lambda_1,\ldots,\Lambda_h}(2^{\nu+1})}\le \sum_{\nu\ge 0} {1\over 2^\nu} \log{1\over\overline\omega_{\Lambda_1,\ldots,\Lambda_h}(2^{\nu+1})}\le \sum_{\nu\ge 0} {1\over 2^\nu} \log{1\over\widetilde\omega_{\Lambda_1,\ldots,\Lambda_h}(2^{\nu+1})}\,,
$$
but the inequalities are generically strict.

\sm Now we can state and prove our holomorphic simultaneous linearization result, Theorem \rf{TeLinSimIntro}, whose proof will be an adaptation to our case of the existing methods introduced by Brjuno [Brj], and Siegel [S1, S2]
(see also P\"oschel [P\"o]). We recall here the statement of Theorem \rf{TeLinSimIntro}.

\sm\thm{TeLinSim}{Theorem} {\sl Let $f_1, \dots, f_h$ be $h\ge 2$ simultaneously formally linearizable germs of biholomorphism of $\C^n$ fixing the origin and such that their linear parts $\Lambda_1,\dots, \Lambda_h$ are simultaneously diagonalizable. If $f_1,\dots, f_h$ satisfy the simultaneous Brjuno condition, then $f_1,\dots f_h$ are holomorphically simultaneously linearizable.}

\sm\proof Without loss of generality we may assume that $\Lambda_1,\dots\Lambda_h$ are diagonal, that is 
$$
\Lambda_k=\diag(\lambda_{k,1},\dots,\lambda_{k,n}).
$$

Since each~$f_k$ is holomorphic in a neighbourhood of the origin, there exists a positive number~$\rho$ such that~$\|f^{(k)}_L\|\le \rho^{|L|}$ for~$|L|\ge 2$. The functional equation
$$
f_k\circ\phe=\phe\circ\Lambda_k,
$$
remains valid under the linear change of coordinates~$f_k(z)\mapsto \sigma f_k(z/\sigma)$,~$\phe(w)\mapsto \sigma\phe(w/\sigma)$ with~$\sigma=\max\{1, \rho^2\}$. Therefore we may assume that 
$$
\forevery{|L|\ge 2}{\|f^{(k)}_L\|\le 1,}
$$
for $k=1,\dots, h$.

Thanks to Theorem \rf{TeLinNotRes}, we may assume that $\phe_{Q,j}=0$ for all $Q\in\cap_{k=1}^h\res_j(\Lambda_k)$, hence for any multi-index $Q\in\N^n\setminus\cap_{k=1}^h\cap_{j=1}^n\res_j(\Lambda_k)$ with $|Q|\ge 2$ we have 
$$
\|\phe_Q\|\le \eps_{Q}^{-1} \sum_{Q_1+\cdots +Q_\nu = Q \atop \nu\ge 2} \|\phe_{Q_1}\|\cdots \|\phe_{Q_\nu}\|,
\tag{eq7bis}
$$
where $\eps_Q =\min_{1\le j \le n}\max_{1\le k\le h}|\Lambda_{k}^Q - \lambda_{k,j}|$.

\sm Now, the proof follows closely the proof of the main Theorem in [P\"o]. For the benefit of the reader, we shall report it here.

\sm Following P\"oschel [P\"o], we inductively define: 
$$
\cases{
\displaystyle\alpha_m= \sum_{m_1+\cdots + m_\nu =m \atop \nu \ge 2} \alpha_{m_1} \cdots \alpha_{m_\nu}, &\hbox{
for}~$m \ge 2$ 
\cr\noalign{\sm}
\alpha_1 = 1,
}
$$
and
$$
\cases{
\displaystyle\delta_Q = \eps_Q^{-1}\max_{Q_1+\cdots + Q_\nu =Q\atop \nu\ge 2} \delta_{Q_1}\cdots\delta_{Q_\nu}, &\hbox{for}~$Q\in\N^n\setminus\cap_{k=1}^h\cap_{j=1}^n\res_j(\Lambda_k)$~\hbox{with}~$|Q|\ge 2$,
\cr\noalign{\sm}
\delta_E= 1, &\hbox{for any}~$E\in\N^n$~\hbox{with}~$|E|=1$.
}$$
Then, by induction, we have
$$
{\|\phe_Q\|\le \alpha_{|Q|}\delta_Q,}
$$
for every $Q\in\N^n\setminus\cap_{k=1}^h\cap_{j=1}^n\res_j(\Lambda_k)$ with~$|Q|\ge 2$. Therefore, to establish 
$$
\sup_Q{1\over |Q|} \log\|\phe_Q\|<+\io,
\tag{eq6}
$$
it suffices to prove analogous estimates for~$\alpha_m$ and~$\delta_Q$.

\sm It is easy to estimate~$\alpha_m$, and we refer to [P\"o] (see also [R4]) for a detailed proof of
$$
\sup_m {1\over m}\log \alpha_m <+\io.
$$

\sm To estimate~$\delta_Q$ we have to take care of small divisors.
First of all, for each multi-index~$Q\not\in\cap_{k=1}^h\cap_{j=1}^n\res_j(\Lambda_k)$ with~$|Q|\ge 2$ we can associate to~$\delta_Q$ a decomposition of the form
$$
\delta_Q= \eps_{L_0}^{-1}\eps_{L_1}^{-1}\cdots\eps_{L_p}^{-1},\tag{eqdelta}
$$
where~$L_0=Q$,~$|Q|>|L_1|\ge\cdots\ge|L_p|\ge2$ and $L_a\not\in\cap_{k=1}^h\cap_{j=1}^n\res_j(\Lambda_k)$ for all $a=1,\dots, p$ and $p\ge 1$. 
In fact, we choose a decomposition~$Q=Q_1+\cdots+Q_\nu$ such that the maximum in the expression of~$\delta_Q$ is achieved; obviously, $Q_a$ does not belong to $\cap_{k=1}^h\cap_{j=1}^n\res_j(\Lambda_k)$ for all $a=1,\dots,\nu$. We can then express~$\delta_Q$ in terms of~$\eps_{Q_j}^{-1}$ and~$\delta_{Q'_j}$  with~$|Q'_j|<|Q_j|$. Carrying on this process, we eventually arrive at a decomposition of the form \rf{eqdelta}. Furthermore, for each multi-index $Q\not\in\cap_{k=1}^h\cap_{j=1}^n\res_j(\Lambda_k)$ with $|Q|\ge 2$, we can choose an index $k_Q$ and an index $i_Q$ so that
$$
\eps_Q = |\Lambda_{k_Q}^Q - \lambda_{k_Q,i_Q}|.
$$ 

For $m\ge 2$ and $1\le j\le n$, we can define
$$
N^j_m(Q)
$$
to be the number of factors~$\eps_{L}^{-1}$ in the expression \rf{eqdelta} of~$\delta_Q$, satisfying
$$
\eps_{L}<\theta\,\omega_{f_1,\dots,f_h}(m),~~\hbox{and}~~i_L = j,
$$ 
where~$\omega_{f_1,\dots,f_h}(m)$ is defined as
$$
\omega_{f_1,\dots,f_h}(m)= \min_{2\le|Q|\le m\atop Q\not\in\cap_{k=1}^h\cap_{j=1}^n\res_j(\Lambda_k)}\eps_Q, 
$$
and~$\theta$ is the positive real number satisfying
$$
4\theta=\min_{1\le p\le n\atop 1\le k\le h }|\lambda_{k,p}| \le 1.
$$
The last inequality can always be satisfied by replacing~$f_k$ by~$f_k^{-1}$ if necessary. Moreover we also have~$\omega_{f_1,\dots,f_h}(m)\le 2$.
Notice that~$\omega_{f_1,\dots,f_h}(m)$ is non-increasing with respect to~$m$ and under our assumptions~$\omega_{f_1,\dots,f_h}(m)$ tends to zero as~$m$ goes to infinity. 

\sm The following is the key estimate, and it descends from Brjuno.

\sm\thm{Le1.1}{Lemma} {\sl For~$m\ge2$,~$1\le j\le n$,~and $Q\not\in\cap_{k=1}^h\cap_{j=1}^n\res_j(\Lambda_k)$, we have
$$
N^j_m(Q)\le \cases{0, &\hbox{if}~~$|Q|\le m$,\cr\noalign{\sm}
					\displaystyle {2|Q|\over m}-1,&\hbox{if}~~$|Q|> m$.}
$$}

\sm The proof of Lemma \rf{Le1.1} can be obtained adapting the proof of Brjuno's lemma contained in the addendum of [P\"o], and we report it here for the sake of completeness.

\sm\proof The proof is done by induction on $|Q|$. Since we fix~$m$ and~$j$ throughout the proof, we shall write~$N$ instead of~$N^j_m$.

For~$|Q|\le m$, 
$$
\eps_Q\ge\omega_{f_1,\dots,f_h}(|Q|)\ge \omega_{f_1,\dots,f_h}(m) > \theta\, \omega_{f_1,\dots,f_h}(m),
$$
hence~$N(Q)=0$.

Assume now that~$|Q|>m$. Then~$2|Q|/m -1 \ge 1$. Write
$$
\delta_Q= \eps_Q^{-1}\delta_{Q_1}\cdots \delta_{Q_\nu},
$$
with
$$
Q=Q_1 + \cdots + Q_\nu, \quad \nu\ge2,\quad\hbox{and}\quad |Q|>|Q_1|\ge \cdots\ge|Q_\nu|;
$$ 
notice that $Q-Q_1$ does not belong to $\cap_{k=1}^h\cap_{j=1}^n\res_j(\Lambda_k)$, otherwise the others $Q_h$'s would be in $\cap_{k=1}^h\cap_{j=1}^n\res_j(\Lambda_k)$.

\no In this decomposition, only $|Q_1|$ can be greater than $M:= \max(|Q|-m, m)$. If this is the case, we can decompose $\delta_{Q_1}$ in the same way, and repeating this step at most $m-1$ times, we obtain the decomposition (where $P_1 = Q_1$)
$$
\delta_Q= \eps_Q^{-1}\eps_{P_1}^{-1}\cdots \eps_{P_\mu}^{-1}\delta_{L_1}\cdots \delta_{L_\nu},
$$
where $\mu\ge 1$, $\nu\ge 2$ and
$$
\eqalign{
&Q > P_1 >\cdots > Q_\mu,\cr
&L_1+\cdots+ L_\nu = Q,\cr
&|P_\mu|> M\ge |L_1|\ge \cdots \ge |L_\nu|. 
}
$$
Here, $Q>L$ means, as in [P\"o], that $Q-L\in \N^n$ is not identically zero. The crucial point is that at most {\sl one} of the $\eps_K$'s can contribute to $N(Q)$, which is the content of the following lemma descending from Siegel.

\sm\thm{Le1.2}{Lemma} {\sl If~$Q > L$, the multi-indices $Q$, $L$ and $Q-L$ are not in $\cap_{k=1}^h\cap_{j=1}^n\res_j(\Lambda_k)$, and
$$
\eps_Q < \theta\, \omega_{f_1,\dots,f_h}(m), \quad \eps_{L}<\theta\, \omega_{f_1,\dots,f_h}(m),\quad i_Q= i_L, 
$$
then~$|Q-L| \ge m$.}

\sm The proof of Lemma \rf{Le1.2} can be obtained adapting the proof of Siegel's lemma contained in the addendum of [P\"o], and we report it here for the sake of completeness.

\sm\proof Thanks to the definition, for all $k\in\{1, \dots,h\}$ we have 
$$
|\Lambda_k^Q-\lambda_{k,i_Q}| \le |\Lambda_{k_Q}^Q-\lambda_{k_Q,i_Q}| \quad\hbox{and}\quad |\Lambda_k^{L}-\lambda_{k,i_L}|\le |\Lambda_{k_L}^L-\lambda_{k_L,i_L}|;
$$
moreover, since we are supposing~$\eps_L <\theta\,\omega_{f_1,\dots, f_h}(m)$, we have
$$
\eqalign{|\Lambda_k^L| &>|\lambda_{k,i_L}|-\theta\,\omega_{f_1,\dots,f_h}(m)\cr
								&\ge 4\theta - 2\theta = 2\theta.
}$$
It follows that
$$
\eqalign{2\theta\,\omega_{f_1,\dots,f_h}(m) &> \eps_Q+ \eps_L\cr
							&\ge |\Lambda_k^Q-\lambda_{k,i_Q}| + |\Lambda_k^L-\lambda_{k,i_L}|\cr
							&\ge |\Lambda_k^Q - \Lambda_k^L|\cr
							&\ge |\Lambda_k^L| \,|\Lambda_k^{Q-L} - 1|\cr
							&\ge 2\theta\,\omega_{f_1,\dots,f_h}(|Q-L|+1),
}$$
and therefore $\omega_{f_1,\dots,f_h}(|Q-L|+1)<\omega_{f_1,\dots,f_h}(m)$, which implies $|Q-L|\ge m$ by the monotonicity of $\omega_{f_1,\dots,f_h}$. \qed

Thanks to the previous result, it follows from the decomposition of $\delta_Q$ that
$$
N(Q) \le 1 + N(L_1) + \cdots + N(L_\nu).
$$
Choose $0\le \rho\le \nu$ such that $|L_\rho|>m \ge |L_{\rho +1}|$. By the induction hypothesis, all terms with $|L|\le m$ vanish, and we obtain
$$\eqalign{
N(Q) 
&\le 1 + N(L_1) + \cdots + N(L_\rho)\cr
&\le 1 + 2{|L_1 +\cdots + L_\rho|\over m}-\rho\cr
&\le\cases{1, &$\rho=0$\cr\noalign{\sm}
		   2{|Q|-m\over m}, &$\rho=1$\cr\noalign{\sm}
		   2{|L_1 +\cdots + L_\rho|\over m} -1, &$\rho\ge 2$,}\cr
&\le 2{|Q|\over m}-1,
}
$$
concluding the proof.\qed

\sm Since the $f_1,\dots, f_h$ satisfy the simultaneous Brjuno condition, there exists a strictly increasing sequence~$\{p_\nu\}_{\nu\ge 0}$ of integers with~$p_0=1$, and such that
$$
\sum_{\nu\ge 0} {1\over p_\nu}\log{1\over\omega_{f_1,\dots,f_h}(p_{\nu +1})} <+\io.\tag{eq8}
$$
We have to estimate
$$
{1\over |Q|}\log\delta_Q = \sum_{j=0}^p {1\over |Q|} \log \eps_{L_j}^{-1}, \quad Q\not\in\cap_{k=1}^h\cap_{j=1}^n\res_j(\Lambda_k).
$$
By Lemma~\rf{Le1.1}, 
$$
\eqalign{{\rm card}\left\{0\le j\le p : \theta\, \omega_{f_1,\dots,f_h}(p_{\nu +1}) \le \eps_{L_j} <\theta\, \omega_{f_1,\dots,f_h}(p_\nu)\right\} 
							&\le {2n|Q|\over p_\nu}} 
$$
for~$\nu\ge 1$. It is also easy to see from the definition of~$\delta_Q$ that the number of factors~$\eps_{L_j}^{-1}$ is bounded by~$2|Q| - 1$, and so the previous inequality holds also for $\nu=0$ when the upper bound is dropped.
Therefore,
$$
\eqalign{{1\over |Q|} \log \delta_Q &\le 2n \sum_{\nu\ge 0} {1\over p_\nu} \log{1\over\theta\,\omega_{f_1,\dots,f_h}(p_{\nu +1})} \cr
	&= 2n\left( \sum_{\nu \ge 0} {1\over p_\nu}\log{1\over\omega_{f_1,\dots,f_h}(p_{\nu +1})} + \log{1\over\theta} \sum_{\nu \ge 0} {1\over p_\nu}\right).}\tag{eq9}
$$
Since~$\omega_{f_1,\dots,f_h}(m)$ tends to zero monotonically as~$m$ goes to infinity, we can choose some~$\overline{m}$ such that~$1>\omega_{f_1,\dots,f_h}(m)$ for all~$m>\overline{m}$, and we obtain
$$
\sum_{\nu\ge\nu_0} {1\over p_\nu} \le {1\over \log (1/\omega_{f_1,\dots,f_h}(\overline{m}))} \sum_{\nu\ge \nu_0} {1\over p_\nu} \log{1\over\omega_{f_1,\dots,f_h}(p_{\nu +1})},
$$
where~$\nu_0$ verifies the inequalities~$p_{\nu_0 -1}\le \overline{m} < p_{\nu_0}$. Thus both series in parentheses in \rf{eq9} converge thanks to \rf{eq8}. Therefore
$$
\sup_Q {1\over |Q|}\log \delta_Q <+\io
$$ 
and this concludes the proof. \qed

\sm Combining Proposition \rf{TeSimFormLinDiag} and Theorem \rf{TeLinSim} we obtain the following equivalence, presented as Theorem \rf{TeSimHolLinDiagIntro} in the introduction.

\sm\thm{TeSimHolLinDiag}{Theorem} {\sl Let $f_1, \dots, f_h$ be $h\ge 2$  formally linearizable germs of biholomorphisms of $\C^n$ fixing the origin, with simultaneously diagonalizable linear parts, and satisfying the simultaneous Brjuno condition. Then $f_1, \dots, f_h$ are simultaneously holomorphically linearizable if and only if they all commute pairwise.}

\sm\proof If $f_1, \dots, f_h$ are simultaneously holomorphically linearizable, since their linear parts are simultaneously diagonalizable, then they all commute and we are done.

On the other hand, if $f_1, \dots, f_h$ commute pairwise then they are simultaneously formally linearizable, by Proposition \rf{TeSimFormLinDiag}, hence the assertion follows from Theorem \rf{TeLinSim}. \qed

We can also deduce the following result.

\sm\thm{TeHolLin}{Corollary} {\sl Let $f_1, \dots, f_h$ be $h\ge 2$ germs of biholomorphism of $\C^n$ fixing the origin and such that the linear part of $f_1$ is diagonalizable, and its eigenvalues have no resonances. If the linear parts of $f_1,\dots, f_h$ are simultaneously diagonalizable, $f_1$ commutes with $f_k$ for $k=2,\dots,h$ and $f_1,\dots, f_h$ satisfy the simultaneous Brjuno condition, then $f_1,\dots, f_h$ are holomorphically simultaneously linearizable.}

\sm\proof It follows from Proposition \rf{LeSimFormLin} that $f_1, \dots,f_h$ are simultaneously formally linearizable, and then the assertion follows from Theorem \rf{TeLinSim}. \qed

\sm We also obtain the following generalization of Theorem $2.5$ of [R1], to which we refer for the definitions of {\it only level $s$ resonances} and {\it simultaneous osculating manifold}. 

\sm\thm{TeGenOsculante}{Proposition} {\sl Let $f_1,\dots,f_h$ be $h\ge 2$ germs of biholomorphism of $\C^n$ fixing the origin with simultaneously diagonalizable linear parts, and satisfying the simultaneous Brjuno condition. Assume that the spectrum of the linear part of $f_1$ has only level $s$ resonances and that $f_1$ commutes with $f_k$ for $k=2,\dots, h$. Then $f_1,\dots, f_h$ are simultaneously holomorphically linearizable if and only if there exists a germ of complex manifold~$M$ at~$O$ of codimension~$s$, invariant under~$f_h$ for each~$h=1, \dots, m$, which is a simultaneous osculating manifold for~$f_1, \dots, f_m$ and such that~$f_1|_M, \dots, f_m|_M$ are simultaneously holomorphically linearizable.}

\sm\proof Proposition \rf{PropOsculanteFormale} implies that $f_1,\dots, f_h$ are simultaneously formally linearizable, hence they are in the hypotheses of Theorem \rf{TeLinSim} and this concludes the proof. \qed

\me\thm{ReRussmann}{Remark} It is possible to prove a simultaneous holomorphic linearization result using the same functional technique of R\"ussmann [R\"u], yielding the following statement.

\sm\no{\sl Given $f_1, \dots, f_h$, $h\ge 2$ pairwise commuting germs of biholomorphism of $\C^n$ fixing the origin and such that their linear parts $\Lambda_1,\dots, \Lambda_h$ are non-resonant and diagonalizable, if 
$$
\sum_{\nu\ge 0} {1\over 2^\nu} \log{1\over\widetilde\omega_{\Lambda_1,\ldots,\Lambda_h} (2^{\nu+1})}<+\io,
\tag{russmann}
$$
where
$$
\widetilde\omega_{\Lambda_1,\ldots,\Lambda_h} (m) = \max_{1\le k\le h} \min_{2\le |Q|\le m\atop 1\le j \le n}|\Lambda_{k}^Q - \lambda_{k,j}|,
$$ 
then $f_1,\dots f_h$ are holomorphically simultaneously linearizable.}

\sm In fact, by the commutation hypothesis, $\Lambda_1,\dots,\Lambda_h$ are simultaneously diagonalizable, so we may assume them to be diagonal. Moreover, since $\Lambda_1,\dots,\Lambda_h$ are non-resonant, each $f_h$ admits a unique formal linearization which is their unique formal simultaneous linearization $\phe$, thanks to the commutation hypothesis. 

We can then perform the argument of Section 6 of [R\"u] for each $f_k$, getting estimates for the convergence of $\phe$ substituting his function $\Omega$ with the function $\omega_{\Lambda_k}(m):=\min_{2\le |Q|\le m\atop 1\le j \le n}|\Lambda_{k}^Q - \lambda_{k,j}|$ for each $k=1,\dots, h$, and we get the assertion because we can estimate the convergence of $\phe$ using the maximum of these $\omega_{\Lambda_1}(m),\dots, \omega_{\Lambda_h}(m)$, that is using \rf{russmann}.

Note that in this case we use one of the generalizations to a simultaneous Brjuno condition we introduced in Remark \rf{RmkSimBrjuno}, which is, as we remarked, stronger than the simultaneous Brjuno condition we introduced. Moreover, to use Theorem 2.1 of [R\"u]  we have to assume $\Lambda_1,\dots,\Lambda_h$ to be non-resonant to have a unique formal linearization, and to be able to apply R\"ussmann's procedure to each germ separately, which is, again, a bit stronger condition than just asking for our germs to be simultaneously formally linearizable. However, at least in principle, it should be possible to prove Theorem \rf{TeLinSim} with functional arguments.

\appendix Equivalence of various Brjuno-type series

\newcount\sectno\sectno=4

\sm In this appendix we prove the equivalence of three Brjuno-type condition/series. In fact in the literature studying of the linearization problem one can find at least the following three series:
$$
B(\omega)= \sum_{\nu\ge 0} {1\over 2^\nu}\log{1\over \omega(2^{\nu +1})},
$$
$$
R(\omega)= \sum_{k\ge 1} {1\over k^2}\log{1\over \omega(k)},
$$
and
$$
\Gamma(\omega)= \sum_{k\ge 1} {1\over k(k+1)}\log{1\over \omega(k)},
$$
where $\omega\colon\N\to\R$ is a non-increasing monotone function, usually containing the information on the small divisors one has to estimate. 
The first one was introduced by Brjuno in [Brj] while the other two can be found for example in [D], [R\"u] and [GM].

R\"ussmann proved that, in dimension $1$, the convergence of $R(\omega)$ is equivalent to the convergence of $B(\omega)$ (see Lemma~8.2 of [R\"u]), and he also proved the following result.

\sm\thm{LeRussmann}{Lemma} (R\"ussmann, 2002 [R\"u]) {\sl Let $\Omega\colon\N\to(0,+\io)$ be a monotone non decreasing function, and 
let $\{s_\nu\}$ be defined by $s_\nu := 2^{q+\nu}$, with $q\in\N$. Then
$$
\sum_{\nu\ge 0} {1\over s_\nu} \log \Omega(s_{\nu+1}) \le \sum_{k\ge 2^{q+1}} {1\over k^2} \log \Omega(k).
$$ 
}

\sm We proved (see Lemma 4.2 of [R4]) that in fact if $R(\omega)<+\io$ then also $B(\omega)<+\io$. In fact, we can prove the following equivalence.

\sm\thm{TeEquivalenza}{Theorem} {\sl Let $\omega\colon\N\to(0,1)$ be a non-increasing monotone function and consider the following series
$$
B(\omega)= \sum_{\nu\ge 0} {1\over 2^\nu}\log{1\over \omega(2^{\nu +1})},\quad
R(\omega)= \sum_{k\ge 1} {1\over k^2}\log{1\over \omega(k)},~\hbox{and}~~
\Gamma(\omega)= \sum_{k\ge 1} {1\over k(k+1)}\log{1\over \omega(k)}.
$$
Then we have the following inequalities:
$$
\Gamma(\omega)\le R(\omega)\le 2\,\Gamma(\omega),\tag{due}
$$
and
$$
\Gamma(\omega) \le {1\over 2} B(\omega)\le 2\, \Gamma(\omega) - \log{1\over\omega(1)}.\tag{uno}
$$
}

\sm\proof The first inequalities are immediate to recover, because for any $k\ge 1$ we have
$$
{1\over k(k+1)} \le {1\over k^2}\le {2\over k(k+1)}.
$$
To prove \rf{uno}, recall that 
$$
\sum_{k=a}^{b-1} {1\over k(k+1)} = {1\over a} - {1\over b}.
$$
Then, for the first inequality in \rf{uno} we have that 
$$
{1\over \omega(2^\nu)} \le \cdots\le {1\over\omega(2^{\nu+1} - 1)} \le{1\over\omega(2^{\nu+1})},  
$$
hence
$$
\eqalign{
\Gamma(\omega)
&=\sum_{\nu\ge0}\sum_{k=2^{\nu}}^{2^{\nu+1}-1} {1\over k(k+1)}\log{1\over \omega(k)}\cr
&\le\sum_{\nu\ge 0}\log{1\over \omega(2^{\nu+1})}\sum_{k=2^{\nu}}^{2^{\nu+1}-1} {1\over k(k+1)}\cr
&={1\over 2} B(\omega).
}
$$
For the other inequality, we use
$$
{1\over \omega(2^{\nu+1})} \le \cdots\le {1\over\omega(2^{\nu+2}- 1)} \le{1\over\omega(2^{\nu+2})};  
$$
hence
$$
\eqalign{
{1\over 2}B(\omega)
&=2\sum_{\nu\ge 0}\log{1\over \omega(2^{\nu+1})}\sum_{k=2^{\nu+1}}^{2^{\nu+2}-1} {1\over k(k+1)}\cr
&\le 2\sum_{\nu\ge 0}\sum_{k=2^{\nu+1}}^{2^{\nu+2}-1} {1\over k(k+1)}\log{1\over \omega(k)}\cr
&=2\,\Gamma(\omega) - \log{1\over \omega(1)}\,,
}
$$
and we are done. \qed

\sm As a corollary we obtain the following equivalence.

\sm\thm{CoEquivalenza}{Corollary} {\sl With the above definitions, we have
$$
B(\omega)< + \io \iff \Gamma(\omega)< + \io \iff R(\omega)< + \io.
$$
}

\sm In particular, the simultaneous Brjuno condition, or any Brjuno condition, can be expressed through the convergence of any of the three series $B$, $R$, and $\Gamma$. 

\vbox{\vskip.75truecm\advance\hsize by 1mm
	\hbox{\centerline{\sezfont References}}
	\vskip.25truecm}\nobreak
{\parindent=40pt

\bib{A1} {\sc Abate, M.:} {\sl Discrete local holomorphic dynamics,} in ``Proceedings of 13th Seminar of Analysis and its Applications, Isfahan, 2003'', Eds. S. Azam et al., University of Isfahan, Iran, 2005, pp. 1--32.

\bib{A2} {\sc Abate, M.:} {\sl Discrete holomorphic local dynamical systems,} to appear in ``Holomorphic Dynamical Systems'', G. Gentili, J. Guenot, G. Patrizio eds., Lectures notes in Math., Springer Verlag, 
Berlin, 2010, arXiv:0903.3289v1.
\bib{Ar} {\sc Arnold, V. I.:} ``Geometrical methods in the theory of ordinary differen\-tial equations'', Springer-Verlag, Berlin, 1988.
\bib{Bi} {\sc Biswas, K.:} {\sl Simultaneous linearization of commuting germs of holomorphic diffeomorphisms}, Preprint, 2009, arXiv:0911.2766v2

\bib{Bra} {\sc Bracci, F.:} {\sl Local dynamics of holomorphic diffeomorphisms,} Boll. UMI (8), 7--B (2004), pp. 609--636.
\bib{Brj} {\sc Brjuno, A. D.:} {\sl Analytic form of differential equations,} Trans. Moscow Math. Soc., {\bf 25} (1971), pp. 131--288; {\bf 26} (1972), pp. 199--239.

\bib{D} {\sc DeLatte, D.}: {\sl Diophantine conditions for the linearization of commuting holomorphic functions,} Disc. Cont. Dyn. Syst. {\bf 3} (1997), pp. 317--332.

\bib{DG} {\sc DeLatte, D., Gramchev, T.}: {\sl Biholomorphic maps with linear parts having Jordan blocks: linearization and resonance type phenomena}, Math. Phys. Electron. J. {\bf 8} (2002), Paper 2, 27 pp.

\bib{\'E} {\sc \'Ecalle, J.}: {\sl Singularit\'es non abordables par la g\'eom\'etrie}, Ann. Inst. Fourier (Grenoble) {\bf 42} (1992), no. 1-2, pp. 73--164. 
\bib{FK} {\sc Fayad, B., Khanin, K.:} {\sl Smooth linearization of commuting circle diffeomorphisms}, Ann. of Math. (2) {\bf 170} (2009), no. 2, pp. 961--980. 
\bib{GM} {\sc Giorgilli, A., Marmi, S.}: {\sl Improved estimates for the convergence radius in the Poincar\'e-Siegel problem}, Discrete and Continuous Dynamical Systems series S {\bf 3}, (2010), pp.601--621.

\bib{GY} {\sc Gramchev, T., Yoshino, M.:} {\sl Rapidly convergent iteration method for simultaneous normal forms of commuting maps},  Math. Z. {\bf 231}  (1999),  no. 4, pp. 745--770. 

\bib{H} {\sc Herman, M.-R.:} {\sl Sur la conjugaison diff\'erentiable des diff\'eomorphismes du cercle \`a des rotations,} Inst. Hautes \'Etudes Sci. Publ. Math. {\bf 49} (1979), pp. 5--233. 

\bib{M} {\sc Moser, J.:} {\sl On commuting circle mappings and simultaneous Diophantine approximations,} Math. Z. {\bf 205} (1990), no. 1, pp. 105--121. 
\bib{P} {\sc P\'erez-Marco, R.}: {\sl Fixed points and circle maps}, Acta Math., {\bf 179}, (1997), pp. 243--294.

\bib{P\"o}{\sc P\"oschel, J.}: {\sl On invariant manifolds of complex analytic mappings near fixed points}, Exp. Math. {\bf 4}, (1986), pp. 97–-109.

\bib{R1} {\sc Raissy, J.}: {\sl Linearization of holomorphic germs with quasi-Brjuno fixed points}, Mathematische Zeitschrift, {\bf 264}, (2010), pp. 881--900.

\bib{R2} {\sc Raissy, J.}: {\sl Simultaneous linearization of holomorphic germs in presence of resonances}, Conform. Geom. Dyn. {\bf 13} (2009), pp. 217--224. 
\bib{R3} {\sc Raissy, J.}: {\bf Geometrical methods in the normalization of germs of biholomorphisms}, Ph.D Thesis, Universit\`a di Pisa (2009).

\bib{R4} {\sc Raissy, J.}: {\sl Brjuno conditions for linearization in presence of resonances}, in  {\bf ``Asymptotics in Dynamics, Geometry and PDE's; Generalized Borel Summation'' vol. I}, O. Costin, F. Fauvet, F. Menous e D. Sauzin editors, {``CRM series'', Pisa, Edizioni Della Normale 2011}, (2011), pp. 201--218.  

\bib{S1} {\sc Siegel, C. L.}: {\sl Iteration of analytic functions}, Ann. of Math. (2) {\bf 43}, (1942), pp. 607-–612.

\bib{S2} {\sc Siegel, C. L.}: {\sl \"Uber die Normalform analytischer Differentialgleichungen in der N\"ahe einer Gleichgewichtsl\"osung}, Nachr. Akad. Wiss. G\"ottingen. Math.-Phys. Kl. (1952), pp. 21-–30.

\bib{R\"u}{\sc R\"ussmann, H.}: {\sl Stability of elliptic fixed points of analytic area-preserving mappings under the Brjuno condition}, Ergodic Theory Dynam. Systems, {\bf 22}, (2002), pp. 1551--1573.

\bib{Y1} {\sc Yoccoz, J.-C.:} {\sl Conjugaison diff\'erentiable des diff\'eomorphismes du cercle dont le nombre de rotation v\'erifie une condition diophantienne}, Ann. Sci. \'Ecole Norm. Sup. {\bf 17} (1984), pp. 333--359. 

\bib{Y2} {\sc Yoccoz, J.-C.:} {\sl Th\'eor\`eme de Siegel, nombres de Bruno et polyn\^omes quadratiques,} Ast\'erisque {\bf 231} (1995), pp. 3--88.
\bib{Y3} {\sc Yoccoz, J.-C.:} {\sl Analytic linearization of circle diffeomorphisms}, in ``Dynamical Systems and Small Divisors (Cetraro, 1998)'', Lecture Notes in Math. 1784, Springer-Verlag, New York, 2002, pp. 125--173. 

\bib{Yo} {\sc Yoshino, M.:} {\sl Diophantine phenomena in commuting vector fields and diffeomorphisms}, Tsukuba J. Math. {\bf 28} (2004), no. 2, pp. 389--399.

}  

\vfill\eject

\no{\titfont Erratum}

\me\no{\sezfont Holomorphic Linearization of commuting Germs}

\sm\no{\sezfont  of Holomorphic Maps }

\bi\no{\bf Jasmin Raissy}

\bi My original paper failed to cite the pioneering
work of Laurent Stolovitch in the paper [S]. In particular, Theorem 1.3 can be seen as a consequence of Theorem 2.1 in [S].

\vbox{\vskip.75truecm\advance\hsize by 1mm
	\hbox{\centerline{\sezfont References}}
	\vskip.25truecm}\nobreak
\parindent=40pt

\bib{S} {\sc Stolovitch, L.:} {\sl Family of intersecting totally real manifolds of $(\C^n,0)$ and CR-singularities}, Preprint, 2005, arXiv:math/0506052v2.

\bye